\documentclass{amsart}

\usepackage[english]{babel}
\usepackage{epsfig}
\usepackage{latexsym, amssymb, amsmath}
\usepackage[only,fivrm]{rawfonts}
\usepackage{textcomp}
\usepackage{amscd}
\usepackage{amssymb}
\usepackage{amsthm}
\usepackage{amsmath}
\usepackage{tikz}
\usepackage{ifthen}
\usetikzlibrary{snakes,arrows,shapes}

\usepackage{xcolor}
\usepackage{color}

\definecolor{darkgreen}{rgb}{0,.4,0.2}
\definecolor{darkagenta}{rgb}{.5,0,.5}
\definecolor{darkred}{rgb}{0.85,0,0}
\definecolor{darkblue}{rgb}{0,0,.6}
\definecolor{lightgray}{gray}{.95}

\newtheorem{theorem}{Theorem}
\theoremstyle{plain}

\newtheorem{lemma}{Lemma}

\newtheorem{proposition}{Proposition}

\newcommand{\eps}{{\varepsilon}}
\newcommand{\R}{{\mathbb R}}

\newcommand{\ts}{\vspace*{.2cm}}   

\newcommand{\D}{\Delta^*}

\newcommand{\AB}{ \left( \begin{smallmatrix} A\\ B \end{smallmatrix} \right)}
\newcommand{\BA}{ \left( \begin{smallmatrix} B\\ A \end{smallmatrix} \right)}
\newcommand{\zz}{ \left( \begin{smallmatrix}\! 0\!\\\! 0\! \end{smallmatrix} \right)}
\newcommand{\zo}{ \left( \begin{smallmatrix}\! 0\!\\\! 1\! \end{smallmatrix} \right)}
\newcommand{\oo}{ \left( \begin{smallmatrix}\! 1\!\\\! 1\! \end{smallmatrix} \right)}
\newcommand{\oz}{ \left( \begin{smallmatrix}\! 1\!\\\! 0\! \end{smallmatrix} \right)}
\newcommand{\thAB}{ \left( \begin{smallmatrix} \theta(A)\!\\ \theta(B)\! \end{smallmatrix} \right)}
\newcommand{\Mth}{ \left( \begin{smallmatrix} 2\,1\\ 2\,3 \end{smallmatrix} \right)}

\newcommand{\pq}{ \left( \begin{smallmatrix} p\\ q \end{smallmatrix} \right)}
\newcommand{\uv}{ \left( \begin{smallmatrix} u^*\\ v^* \end{smallmatrix} \right)}
\newcommand{\T}{{\rm T}}
\newcommand{\Bal}[2]{\Big|\begin{matrix} #1\\[-.1cm] #2 \end{matrix} \Big|}
\newcommand{\Col}[2]{\Big( \begin{matrix} #1\\[-.1cm] #2 \end{matrix} \Big) }
\newcommand{\tht}{\theta}
\newcommand{\thtt}{\theta\times\!\theta\,}
\newcommand{\hA}{\,\,\hat{\!\!A}}

\newcommand{\hC}{\hat{C}}
\newcommand{\hD}{\hat{D}}
\newcommand{\dd}{\delta}
\newcommand{\ddp}{\dot{\delta}}
\newcommand{\ha}{\dot{h}}
\newcommand{\hb}{\ddot{h}}

\newcommand{\skewzero}{
\begin{tikzpicture}[scale=0.7,rounded corners=2]
\path (0.25,0) coordinate (HS);   \path (0,0.4) coordinate (VS);
\path (VS)++(HS) coordinate  (PVH);   \path (PVH)++(VS) coordinate  (PVHV);
\path (PVH)++(HS) coordinate  (PHVH);   \path (PHVH)++(VS) coordinate  (PHVHV);
\draw (0,0)--(VS)--(PVH)--(PVHV);   \draw (HS)--(PVH)--(PHVH)--(PHVHV);
\path (0.125,0.2) coordinate  (PL);   \path (PL)++(PVH) coordinate  (PLVH);
\node at (PL) { $0$};  \node at (PLVH) {$0$};
\end{tikzpicture}}

\newcommand{\skewone}{
\begin{tikzpicture}[scale=0.7,rounded corners=2]
\path (0.25,0) coordinate (HS);   \path (0,0.4) coordinate (VS);
\path (VS)++(HS) coordinate  (PVH);   \path (PVH)++(VS) coordinate  (PVHV);
\path (PVH)++(HS) coordinate  (PHVH);   \path (PHVH)++(VS) coordinate  (PHVHV);
\draw (0,0)--(VS)--(PVH)--(PVHV);   \draw (HS)--(PVH)--(PHVH)--(PHVHV);
\path (0.125,0.2) coordinate  (PL);   \path (PL)++(PVH) coordinate  (PLVH);
\node at (PL) { $1$};  \node at (PLVH) {$1$};
\end{tikzpicture}}

\newcommand{\Ahat}{
\begin{tikzpicture}[scale=1.1,rounded corners=2]
\path (0.25,0) coordinate (HS);   \path (0,0.4) coordinate (VS); \path (VS)++(VS) coordinate (VS2); \path (HS)++(HS) coordinate (HS2);
\path (VS)++(HS2) coordinate  (PVH2);   \path (PVH2)++(VS) coordinate  (PVHV2);
  \path (PVH2)++(HS) coordinate  (PHVH2);   \path (PHVH2)++(VS) coordinate  (PHVHV2);
\draw (0,0)--(VS2);  \draw (HS2)--(PVH2)--(PHVH2)--(PHVHV2);
\path (0.125,0.2) coordinate  (PL);   \path (PL)++(VS) coordinate  (PLV);
\node at (PL) { $\;\;1\, 1$};  \node at (PLV) {$\;\;\;\;\,0\, 1\, 0$};
\end{tikzpicture}}

\newcommand{\Bhat}{
\begin{tikzpicture}[scale=1.1,rounded corners=2]
\path (0,0) coordinate (startB); 
\path (0.25,0) coordinate (HS);   \path (0,0.37) coordinate (VS); \path (VS)++(VS) coordinate (VS2);
\path (HS)++(HS) coordinate (HS2); \path (HS2)++(HS) coordinate (HS3);
\path (startB) coordinate  (P0);   \path (P0)++(VS) coordinate (P1); \path (P1)++(HS) coordinate (P2); \path (P2)++(VS) coordinate (P3);
\draw (startB)++(HS3)--++(VS2);  \draw (P0)--(P1)--(P2)--(P3);
\path (startB)++(0.125,0.2) coordinate  (PL);   \path (PL)++(VS)++(HS) coordinate  (PLV);
\node at (PL) {$\;\;\;\;\,0\, 1\, 0$ };  \node at (PLV) {$\;\;1\, 1$};
\end{tikzpicture}}

\colorlet{Color00}{red!70!black} \colorlet{Color01}{green!29!white}
\colorlet{Color10}{yellow!55!white} \colorlet{Color11}{blue!70!black}
\colorlet{Color00d}{red!20!white} \colorlet{Color01d}{green!12!white}
\colorlet{Color10d}{yellow!14!white} \colorlet{Color11d}{blue!40!white}

\def\EL{0.0909090909} 
\def\ELL{0.008264462810} 

\begin{document}

\title{ Diagonals in  2D-tilings and coincidence densities of substitutions}

\author{F.~Michel Dekking}


\begin{abstract}
We study the diagonals of two-dimensional tilings generated by direct product substitutions. The properties of these diagonals are primarily determined by the eigenvalues of the substitution matrix, but also the order of the letters in the substitution plays a role. We show that the diagonals may fail to be uniformly recurrent, and that the frequencies of letters on the diagonal may not exist. We also highlight the connection with the density of coincidences and overlap distributions.

\smallskip


\medskip


\end{abstract}

\maketitle

\section{Introduction}\label{sec: intro}

  The two-dimensional tilings generated by direct product substitutions have a very simple structure, but we will show that the \emph{diagonals} of these tilings are not that simple. Throughout the paper $\tht$ denotes a substitution (morphism) on $\{0,1\}$ that admits two fixed points:
 we will require that the first letter of $\tht(0)$ is $0$ and the first letter of $\tht(1)$ is $1$.

 The \emph{(direct) product substitution} $\thtt$ of $\tht$ with itself is defined by
 $$\thtt(i,j)=\tht(i)\times\tht(j), \quad i,j\in\{0,1\},$$
 where for two words $v$ and $w$ of length $m$ and $n$ we define the product word as the matrix
 $$v\times w= \big((v_k,w_\ell)\big)_{1\le k\le m, 1\le \ell\le n}.$$
 For example, when $\tht$ is the Thue-Morse substitution given by $\tht(0)=0\,1$, and $\tht(1)=1\,0$ then,
 $$\thtt((0,0))= \left( \begin{matrix} (0,1)\,(1,1)\\ (0,0)\,(1,0) \end{matrix} \right),\quad{\rm and} \quad
  \thtt( (0,1))= \left( \begin{matrix} (0,0)\,(1,0)\\ (0,1)\,(1,1) \end{matrix} \right).$$

\noindent It is visually more convenient to represent this with $1\times 1$ colored tiles:\\[-.6cm]
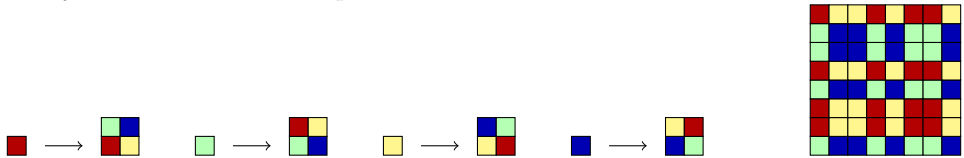
\begin{figure}[h!]\label{fig:Morseprod}
\begin{tikzpicture}[scale=.25,rounded corners=0]
\def\brickzerozero(#1){\path (#1) coordinate (P0);  \draw [fill=Color00, draw=black]  (P0) rectangle +(1,1)}
\def\brickzeroone(#1){\path (#1) coordinate (P0);  \draw [fill=Color01, draw=black]  (P0) rectangle +(1,1)}
\def\brickonezero(#1){\path (#1) coordinate (P0); \draw [fill=Color10, draw=black]  (P0) rectangle +(1,1)}
\def\brickoneone(#1){\path (#1) coordinate (P0);  \draw [fill=Color11, draw=black]  (P0) rectangle +(1,1)}
\def\sigma(#1,#2){\path (#2) coordinate (P0);
\path (P0)++(0,0) coordinate  (P000); \path (P0)++(1,0) coordinate  (P001);
\path (P0)++(0,1) coordinate  (P002); \path (P0)++(1,1) coordinate  (P003);
\ifthenelse{#1=00}{\brickzerozero(P000);\brickonezero(P001);\brickzeroone(P002);\brickoneone(P003)};{}; 
\ifthenelse{#1=01}{\brickzeroone(P000);\brickoneone(P001);\brickzerozero(P002);\brickonezero(P003)};{}; 
\ifthenelse{#1=10}{\brickonezero(P000);\brickzerozero(P001);\brickoneone(P002);\brickzeroone(P003)};{}; 
\ifthenelse{#1=11}{\brickoneone(P000);\brickzeroone(P001);\brickonezero(P002);\brickzerozero(P003)};{}; 
}; 
\path (0,0) coordinate (B00); \path (5,0) coordinate (SB00);
\path (10,0) coordinate (B01); \path (15,0) coordinate (SB01);
\path (20,0) coordinate (B10); \path (25,0) coordinate (SB10);
\path (30,0) coordinate (B11); \path (35,0) coordinate (SB11);
\brickzerozero(B00); \draw [->] (2,0.5) -- (4,0.5); \sigma(00,SB00);
\brickzeroone(B01); \draw [->] (12,0.5) -- (14,0.5); \sigma(01,SB01);
\brickonezero(B10); \draw [->] (22,0.5) -- (24,0.5); \sigma(10,SB10);
\brickoneone(B11); \draw [->] (32,0.5) -- (34,0.5); \sigma(11,SB11);
\end{tikzpicture}\hspace*{-1.2cm}
\begin{tikzpicture}[scale=.25,rounded corners=0]
\path (2,0) coordinate (TH1); \path (0,2) coordinate (TV1); \path (4,0) coordinate (TH2); \path (0,4) coordinate (TV2);
\path (8,0) coordinate (TH3); \path (0,8) coordinate (TV3);
\def\brickzerozero(#1){\path (#1) coordinate (P0);  \draw [fill=Color00, draw=black]  (P0) rectangle +(1,1)}
\def\brickzeroone(#1){\path (#1) coordinate (P0);  \draw [fill=Color01, draw=black]  (P0) rectangle +(1,1)}
\def\brickonezero(#1){\path (#1) coordinate (P0); \draw [fill=Color10, draw=black]  (P0) rectangle +(1,1)}
\def\brickoneone(#1){\path (#1) coordinate (P0); \draw [fill=Color11, draw=black]  (P0) rectangle +(1,1)}
\path (0,0) coordinate (start); \path (2,0) coordinate (next1); \path (4,0) coordinate (next2);  \path (6,0) coordinate (next3);
\def\sigma(#1,#2){\path (#2) coordinate (P0);
\path (P0)++(0,0) coordinate  (P000); \path (P0)++(1,0) coordinate  (P001);
\path (P0)++(0,1) coordinate  (P002); \path (P0)++(1,1) coordinate  (P003);
\ifthenelse{#1=00}{\brickzerozero(P000);\brickonezero(P001);\brickzeroone(P002);\brickoneone(P003)};{}; 
\ifthenelse{#1=01}{\brickzeroone(P000);\brickoneone(P001);\brickzerozero(P002);\brickonezero(P003)};{}; 
\ifthenelse{#1=10}{\brickonezero(P000);\brickzerozero(P001);\brickoneone(P002);\brickzeroone(P003)};{}; 
\ifthenelse{#1=11}{\brickoneone(P000);\brickzeroone(P001);\brickonezero(P002);\brickzerozero(P003)};{}; 
}; 
\def\sigmatwo(#1,#2){\path (#2) coordinate (Q0);
\path (Q0)++(0,0) coordinate (Q000); \path (Q000)++(TH1) coordinate (Q001); \path (Q000)++(TV1) coordinate (Q002); \path (Q002)++(TH1)
coordinate (Q003);
\path (Q0)++(0,0) coordinate (Q010); \path (Q010)++(TH1) coordinate (Q011); \path (Q010)++(TV1) coordinate (Q012); \path (Q012)++(TH1)
coordinate (Q013);
\path (Q0)++(0,0) coordinate (Q100); \path (Q100)++(TH1) coordinate (Q101); \path (Q100)++(TV1) coordinate (Q102); \path (Q102)++(TH1)
coordinate (Q103);
\path (Q0)++(0,0) coordinate (Q110); \path (Q110)++(TH1) coordinate (Q111); \path (Q110)++(TV1) coordinate (Q112); \path (Q112)++(TH1)
coordinate (Q113);
\ifthenelse{#1=00}
{\sigma(00,Q000);\sigma(10,Q001);\sigma(01,Q002);\sigma(11,Q003)} {};   
\ifthenelse{#1=01}
{\sigma(01,Q010);\sigma(11,Q011);\sigma(00,Q012);\sigma(10,Q013)} {};   
\ifthenelse{#1=10}
{\sigma(10,Q100);\sigma(00,Q101);\sigma(11,Q102);\sigma(01,Q103)} {};   
\ifthenelse{#1=11}
{\sigma(11,Q110);\sigma(01,Q111);\sigma(10,Q112);\sigma(00,Q113)} {};   
}; 
\def\sigmathree(#1,#2){\path (#2) coordinate (R0);
\path (R0)++(0,0) coordinate (R000); \path (R000)++(TH2) coordinate (R001); \path (R000)++(TV2) coordinate (R002); \path (R002)++(TH2)
coordinate (R003);
\path (R0)++(0,0) coordinate (R010); \path (R010)++(TH2) coordinate (R011); \path (R010)++(TV2) coordinate (R012); \path (R012)++(TH2)
coordinate (R013);
\path (R0)++(0,0) coordinate (R100); \path (R100)++(TH2) coordinate (R101); \path (R100)++(TV2) coordinate (R102); \path (R102)++(TH2)
coordinate (R103);
\path (R0)++(0,0) coordinate (R110); \path (R110)++(TH2) coordinate (R111); \path (R110)++(TV2) coordinate (R112); \path (R112)++(TH2)
coordinate (R113);
\ifthenelse{#1=00}
{\sigmatwo(00,R000);\sigmatwo(10,R001);\sigmatwo(01,R002);\sigmatwo(11,R003)} {};   
\ifthenelse{#1=01}
{\sigmatwo(01,R010);\sigmatwo(11,R011);\sigmatwo(00,R012);\sigmatwo(10,R013)} {};   
\ifthenelse{#1=10}
{\sigmatwo(10,R100);\sigmatwo(00,R101);\sigmatwo(11,R102);\sigmatwo(01,R103)} {};   
\ifthenelse{#1=11}
{\sigmatwo(11,R110);\sigmatwo(01,R111);\sigmatwo(10,R112);\sigmatwo(00,R113)} {};   
}; 
\path (10,0) coordinate (start); \path (start)++(8,8) coordinate (upright);
\sigmathree(01,start); 
\end{tikzpicture}
\caption{The Thue Morse direct product substitution, and an order 3 iterate.}
\end{figure}

\medskip

 Let $u^*=\tht^\infty(0)$ and $v^*=\tht^\infty(1)$ be the two fixed points of $\tht$. Obviously the 2D-substitution $\thtt$ has four fixed
 points,  $u^*\times u^*$, $v^*\times v^*$, $u^*\times v^*$, and $v^*\times u^*$. The first two 2D infinite words have trivial
 diagonals: these are just a doubling of $u^*$, respectively $v^*$. However, even though the diagonal of $u^*\times v^*$ is  also generated by a substitution of length $q$ when $\tht$ is a substitution of constant length $q$, this diagonal is not always so regular as in the $u^*\times u^*$ case.

 From now on, we will often denote the letters of the alphabet $\{0,1\}\times\{0,1\}$ as column vectors, for easier reading. It is easy to see that the  diagonal $\D=\D(0,1)$ of the $u^*\times v^*$ fixed point of the Thue-Morse direct product substitution is generated by the substitution
 $$\Col{0}{1}\mapsto\Col{0}{1}\Col{1}{0}\quad{\rm and} \quad \Col{1}{0}\mapsto\Col{1}{0}\Col{0}{1};
 \quad{\rm in\; tiles:} \quad \begin{tikzpicture}[scale=.2,rounded corners=0]
\def\brickzeroone(#1){\path (#1) coordinate (P0);  \draw [fill=Color01, draw=black]  (P0) rectangle +(1,1)}
\def\brickonezero(#1){\path (#1) coordinate (P0); \draw [fill=Color10, draw=black]  (P0) rectangle +(1,1)}
\def\sigma(#1,#2){\path (#2) coordinate (P0);
\path (P0)++(0,0) coordinate  (P000); \path (P0)++(1,1) coordinate  (P003);
\ifthenelse{#1=01}{\brickzeroone(P000);\brickonezero(P003)};{}; 
\ifthenelse{#1=10}{\brickonezero(P000);\brickzeroone(P003)};{}; 
}; 
\path (15,0) coordinate (B01); \path (20,0) coordinate (SB01);
\path (28,0) coordinate (B10); \path (33,0) coordinate (SB10);
\brickzeroone(B01); \draw [->] (17,0.5) -- (19,0.5); \sigma(01,SB01);
\brickonezero(B10); \draw [->] (30,0.5) -- (32,0.5); \sigma(10,SB10);
\end{tikzpicture}.$$

 We see that the successive first coordinates of the diagonal will form $u^*$ and the successive second coordinates of the diagonal will form  $v^*$. This will always be the case as soon as $\tht$ is \emph{continuous}, which means that ${\tht(1)}=\overline{\tht(0)}$, where
 we use the  mirror substitution  given by $\overline{0}=1$ and $\overline{1}=0$.
 When $\tht$ is not continuous, it is called \emph{discrete}. The terminology follows \cite{CK1971}. As an example of a discrete substitution
 take
 $$\tht(0)= 0\,1\,1 \quad \tht(1)= 1\,0\,1.$$
The substitution that generates the diagonal is then given by
$$\Col{0}{0}\mapsto\Col{0}{0}\Col{1}{1}\Col{1}{1},\; \Col{0}{1}\mapsto\Col{0}{1}\Col{1}{0}\Col{1}{1},\;
  \Col{1}{0}\mapsto\Col{1}{0}\Col{0}{1}\Col{1}{1},\; \Col{1}{1}\mapsto\Col{1}{1}\Col{0}{0}\Col{1}{1}.$$
It follows that for discrete substitutions the diagonal is less regular: the sequence $\D$ is still recurrent, but no longer uniformly
recurrent.
This is caused by the occurrence of longer and longer words consisting solely of the letters  $\zz$ and $\oo$. This also leads to the property that $\zo$ and $\oz$ have frequency 0 and $\zz$ and $\oo$ have a frequency equal to the original frequencies of $0$ and $1$ in the fixed point $u^*$
 of $\tht$ (cf.~Lemma~\ref{lem:subseq}). The technical explanation of this phenomenon is that for discrete substitutions the substitution that generates the diagonal is no longer primitive, but reducible.

\medskip

In Section~\ref{sec: balance} we consider the diagonal problem for substitutions of \emph{non-constant} length. It turns out that in some cases a description of the diagonal is possible with the aid of the notion of \emph{balanced blocks}, which has been used in the literature to compute coincidence densities. In Section~\ref{sec: weights} we extend the notion of balanced blocks, so that it is applicable to a larger class of non-constant substitutions. In Section~\ref{sec: comp} we show how to compute coincidence densities, and give an example where the coincidence density (and so the frequencies of letters) does not exist. In Section~\ref{sec: selfsim} we present a more positive result: using the same  substitution to generate a self-similar tiling, the frequencies of letters on the diagonal do exist. In fact this is shown to be generally true for substitutions with rational eigenvalues. In Section~\ref{sec: outside} we discuss the case of two eigenvalues with modulus larger than 1.
In Section~\ref{sec: overlap} we discuss the connection with the notion of overlap distribution from statistical physics.

\section{Balanced blocks and the density of coincidences}\label{sec: balance}

Non-constant length substitutions can also generate tilings (\cite{Priebe-primer}). The usual example is the Fibonacci substitution
$0\mapsto 0\,1$, $1\mapsto 0$. But since this substitution does not have two fixed points, we will in the following rather consider  its nephew given
by $$\tht:\quad 0\mapsto 0\,0\,1,\quad 1\mapsto 1\,0.$$
 If we want to generate the diagonal $\D$ of the corresponding $u^*\times v^*$ we now have the problem that the
occurrences of $\tht(0)$ and $\tht(1)$ in $u^*$ and $v^*$ are no longer synchronized.
The solution to this problem stems from \cite{PierreM}, and is known as the \emph{Method of balanced blocks}.

Let $N_i(A)$ denote the number of times letter $i$ occurs in the word $A$.
A  `double' word $\AB$ is called a \emph{balanced} block if $$N_0(A)=N_0(B)\;{\rm and}\; N_1(A)=N_1(B).$$
The idea is to cut up the `double' sequence $\D=\uv$ into irreducible balanced blocks:
\begin{equation*}
  \begin{split}
    u^*=\;&  |0\,0\,1|\,0|\,0\,1|\,1\,0|\,0|\,0|\,1|\,0\,0\,1|\,1\,0|\dots \\[-.21cm]
    v^*=\;&  |1\,0\,0|\,0|\,1\,0|\,0\,1|\,0|\,0|\,1|\,1\,0\,0|\,0\,1|\,\dots
  \end{split}
\end{equation*}
The key observation is that for any balanced block $\AB$ the block $\thAB$ is obviously again balanced.
This makes it possible to define a substitution $\hat{\theta}$ on  blocks  $\hA=\AB$ by
$$ \hat{\theta}(\hA)=\hA_1\dots\hA_m,$$
where $\hA_1\dots\hA_m$ is the decomposition of $\hA$ into irreducible balanced blocks.\\
 For example, taking the first balanced block in the decomposition above:
\begin{equation*}
{ \hat{\theta}:}\hspace*{-5cm}\begin{split}
 |0\,0\,1|\;\,_{\mapsto}&\; |0\,0\,1|\,|\,0|\,|\,0\,1|\,|\,1\,0|\\[-.21cm]
 |1\,0\,0|\;\quad          &\;        |1\,0\,0|\,|\,0|\,|\,1\,0|\,|\,0\,1|.
\end{split}
\end{equation*}

Continuing in this way, and identifying the blocks $\AB$ and $\BA$, one obtains $\hat{\theta}$ defined on the alphabet of (identified) irreducible balanced blocks
$$\hat{I}=\Big\{\Bal{0\,0\,1}{1\,0\,0},\, \Bal{0\,1}{1\,0},\,\Bal{0}{0},\,\Bal{1}{1} \Big\}.$$
As in \cite{PierreM} (see also Lemma~\ref{lem:subseq}) it then follows that $\Delta^*$ behaves as in the case of discrete substitutions of
constant length: the symbols  $\hat{0}:=\zz$ and $\hat{1}:=\oo$ dominate.

Actually \cite{PierreM} focusses on the so called \emph{density of coincidences}, defined by
$$\delta_c:=\lim_{n\rightarrow\infty}\,\frac{1}{n}{\rm Card} \{1\le k\le n: u^*_k=v^*_k\}.$$
Since counting coincidences is equivalent to counting occurrences of  $\hat{0}$ and  $\hat{1}$, the result is that $\delta_c=1$ for this substitution.
This plays an important role in the determination of the spectra of dynamical systems generated by the substitutions, see e.g.,
\cite{PierreM}, and \cite{HollSolom}. Crucial there is that  $\tht$ is a substitution of Pisot type, i.e., the second eigenvalue of $M_\tht$ lies in the unit circle. In \cite{Barge-Dia} it is shown that this implies that the number of irreducible balanced blocks is finite.

\medskip

We return to our original question: what is the structure of the diagonal of a direct product substitution?
As an example of a non-Pisot we consider a substitution of Salem type, i.e., the second eigenvalue lies on the unit circle. We take $$\theta:\;0\mapsto 0\,1\,0,\; 1\mapsto 1\,1\,0\,1\,0.$$
Here the `double' fixed point sequence $\D=\uv$ starts with the irreducible balanced block $\Bal{0\,1\,0\,1\,1}{1\,1\,0\,1\,0}$.
Then $\hat{\theta}\big(\Bal{0\,1\,0\,1\,1}{1\,1\,0\,1\,0}\big)$ contains the irreducible balanced block
$\Bal{0\,1\,0\,0\,1\,0\,1\,1}{1\,1\,0\,1\,0\,0\,1\,0}$,
and the  $\hat{\theta}$-image of this block contains a irreducible balanced block of length 27.
One can show that arbitrary long irreducible balanced blocks will occur. In particular we do have a problem to  compute $\delta_c$ for this
$\theta$.

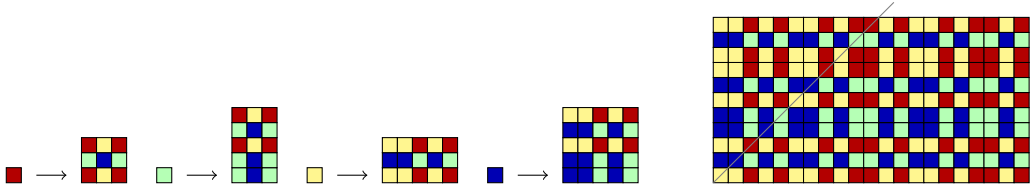
\begin{figure}[t!]\label{fig:weak-mix}
\begin{tikzpicture}[scale=0.2,rounded corners=0]
\path (3,0) coordinate (Th01); \path (5,0) coordinate (Th11); \path (0,3) coordinate (Tv01); \path (0,5) coordinate (Tv11);
\def\brickzerozero(#1){\path (#1) coordinate (P0);  \draw [fill=Color00, draw=black] (P0) rectangle +(1,1)}
\def\brickzeroone(#1){\path (#1) coordinate (P0);  \draw [fill=Color01, draw=black]   (P0) rectangle +(1,1)}
\def\brickonezero(#1){\path (#1) coordinate (P0); \draw [fill=Color10, draw=black]    (P0) rectangle +(1,1)}
\def\brickoneone(#1){\path (#1) coordinate (P0);  \draw [fill=Color11, draw=black]   (P0) rectangle +(1,1)}
\def\sigma(#1,#2){\path (#2) coordinate (P0);
\path (P0)++(0,0) coordinate (P000); \path (P0)++(1,0) coordinate (P001); \path (P0)++(2,0) coordinate (P002);
\path (P0)++(0,1) coordinate (P003); \path (P0)++(1,1) coordinate (P004); \path (P0)++(2,1) coordinate (P005);
\path (P0)++(0,2) coordinate (P006); \path (P0)++(1,2) coordinate (P007); \path (P0)++(2,2) coordinate (P008);
\ifthenelse{#1=00}{\brickzerozero(P000);\brickonezero(P001);\brickzerozero(P002);\brickzeroone(P003);
\brickoneone(P004);\brickzeroone(P005);\brickzerozero(P006);\brickonezero(P007);\brickzerozero(P008)}{}; 
\path (P0)++(0,0) coordinate  (P010); \path (P0)++(1,0) coordinate  (P011); \path (P0)++(2,0) coordinate  (P012);
\path (P0)++(0,1) coordinate  (P013); \path (P0)++(1,1) coordinate  (P014); \path (P0)++(2,1) coordinate  (P015);
\path (P0)++(0,2) coordinate  (P016); \path (P0)++(1,2) coordinate  (P017); \path (P0)++(2,2) coordinate  (P018);
\path (P0)++(0,3) coordinate  (P019); \path (P0)++(1,3) coordinate  (P0110); \path (P0)++(2,3) coordinate  (P0111);
\path (P0)++(0,4) coordinate  (P0112); \path (P0)++(1,4) coordinate  (P0113); \path (P0)++(2,4) coordinate  (P0114);
\ifthenelse{#1=01}
{\brickzeroone(P010);\brickoneone(P011);\brickzeroone(P012);%
\brickzeroone(P013);\brickoneone(P014);\brickzeroone(P015);%
\brickzerozero(P016);\brickonezero(P017);\brickzerozero(P018);%
\brickzeroone(P019);\brickoneone(P0110);\brickzeroone(P0111);%
\brickzerozero(P0112);\brickonezero(P0113);\brickzerozero(P0114);}{}; 
\path (P0)++(0,0) coordinate  (P100); \path (P0)++(1,0) coordinate  (P101); \path (P0)++(2,0) coordinate  (P102);
\path (P0)++(3,0) coordinate  (P103); \path (P0)++(4,0) coordinate  (P104); \path (P0)++(0,1) coordinate  (P105);
\path (P0)++(1,1) coordinate  (P106); \path (P0)++(2,1) coordinate  (P107); \path (P0)++(3,1) coordinate  (P108);
\path (P0)++(4,1) coordinate  (P109); \path (P0)++(0,2) coordinate  (P1010); \path (P0)++(1,2) coordinate  (P1011);
\path (P0)++(2,2) coordinate  (P1012); \path (P0)++(3,2) coordinate  (P1013); \path (P0)++(4,2) coordinate  (P1014);
\ifthenelse{#1=10}
{\brickonezero(P100);\brickonezero(P101);\brickzerozero(P102);\brickonezero(P103);\brickzerozero(P104);%
\brickoneone(P105);\brickoneone(P106);\brickzeroone(P107);\brickoneone(P108);\brickzeroone(P109);%
\brickonezero(P1010);\brickonezero(P1011);\brickzerozero(P1012);\brickonezero(P1013);\brickzerozero(P1014);}{}; 
\path (P0)++(0,0) coordinate (P110); \path (P0)++(1,0) coordinate (P111); \path (P0)++(2,0) coordinate (P112); \path (P0)++(3,0) coordinate
(P113); \path (P0)++(4,0) coordinate (P114);%
\path (P0)++(0,1) coordinate (P115); \path (P0)++(1,1) coordinate (P116); \path (P0)++(2,1) coordinate (P117); \path (P0)++(3,1) coordinate
(P118);\path (P0)++(4,1) coordinate (P119);%
\path (P0)++(0,2) coordinate (P1110); \path (P0)++(1,2) coordinate (P1111);\path (P0)++(2,2) coordinate (P1112); \path (P0)++(3,2) coordinate
(P1113);\path (P0)++(4,2) coordinate (P1114);%
\path (P0)++(0,3) coordinate (P1115); \path (P0)++(1,3) coordinate (P1116);\path (P0)++(2,3) coordinate (P1117); \path (P0)++(3,3) coordinate
(P1118);\path (P0)++(4,3) coordinate (P1119);%
\path (P0)++(0,4) coordinate (P1120); \path (P0)++(1,4) coordinate (P1121);\path (P0)++(2,4) coordinate (P1122); \path (P0)++(3,4) coordinate
(P1123);\path (P0)++(4,4) coordinate (P1124);
\ifthenelse{#1=11}
{\brickoneone(P110);\brickoneone(P111);\brickzeroone(P112);\brickoneone(P113);\brickzeroone(P114);%
\brickoneone(P115);\brickoneone(P116);\brickzeroone(P117);\brickoneone(P118);\brickzeroone(P119);%
\brickonezero(P1110);\brickonezero(P1111);\brickzerozero(P1112);\brickonezero(P1113);\brickzerozero(P1114);%
\brickoneone(P1115);\brickoneone(P1116);\brickzeroone(P1117);\brickoneone(P1118);\brickzeroone(P1119);%
\brickonezero(P1120);\brickonezero(P1121);\brickzerozero(P1122);\brickonezero(P1123);\brickzerozero(P1124);}{}; 
}; 
\path (0,0) coordinate (B00); \path (5,0) coordinate (SB00);
\path (10,0) coordinate (B01); \path (15,0) coordinate (SB01);
\path (20,0) coordinate (B10); \path (25,0) coordinate (SB10);
\path (32,0) coordinate (B11); \path (37,0) coordinate (SB11);
\brickzerozero(B00); \draw [->] (2,0.5) -- (4,0.5); \sigma(00,SB00);
\brickzeroone(B01); \draw [->] (12,0.5) -- (14,0.5); \sigma(01,SB01);
\brickonezero(B10); \draw [->] (22,0.5) -- (24,0.5); \sigma(10,SB10);
\brickoneone(B11); \draw [->] (34,0.5) -- (36,0.5); \sigma(11,SB11);
\def\sigmatwo(#1,#2){\path (#2) coordinate (Q0);
\path (Q0)++(0,0) coordinate  (Q010); \path (Q010)++(Th01) coordinate  (Q011); \path (Q011)++(Th11) coordinate  (Q012);
\path (Q010)++(Tv11) coordinate  (Q013); \path (Q013)++(Th01) coordinate  (Q014); \path (Q014)++(Th11) coordinate  (Q015);
\path (Q013)++(Tv11) coordinate  (Q016); \path (Q016)++(Th01) coordinate  (Q017); \path (Q017)++(Th11) coordinate  (Q018);
\path (Q016)++(Tv01) coordinate  (Q019); \path (Q019)++(Th01) coordinate  (Q0110); \path (Q0110)++(Th11) coordinate  (Q0111);
\path (Q019)++(Tv11) coordinate  (Q0112); \path (Q0112)++(Th01) coordinate  (Q0113); \path (Q0113)++(Th11) coordinate (Q0114);
\ifthenelse{#1=01}
{\sigma(01,Q010);\sigma(11,Q011);\sigma(01,Q012);%
\sigma(01,Q013);\sigma(11,Q014);\sigma(01,Q015);%
\sigma(00,Q016);\sigma(10,Q017);\sigma(00,Q018);%
\sigma(01,Q019);\sigma(11,Q0110);\sigma(01,Q0111);%
\sigma(00,Q0112);\sigma(10,Q0113);\sigma(00,Q0114);}{};                                            
\path (Q0)++(0,0) coordinate (Q100);\path (Q100)++(Th11) coordinate (Q101);\path (Q101)++(Th11) coordinate (Q102);\path (Q102)++(Th01) coordinate (Q103); \path (Q103)++(Th11) coordinate  (Q104);%
\path (Q0)++(Tv01) coordinate (Q105);\path (Q105)++(Th11) coordinate (Q106); \path (Q106)++(Th11) coordinate (Q107);\path (Q107)++(Th01) coordinate (Q108);\path (Q108)++(Th11) coordinate  (Q109);%
\path (Q105)++(Tv11) coordinate (Q1010); \path (Q1010)++(Th11) coordinate (Q1011);\path (Q1011)++(Th11) coordinate (Q1012);\path (Q1012)++(Th01) coordinate (Q1013); \path (Q1013)++(Th11) coordinate  (Q1014);
\ifthenelse{#1=10}
{\sigma(10,Q100);\sigma(10,Q101);\sigma(00,Q102);\sigma(10,Q103);\sigma(00,Q104);%
\sigma(11,Q105);\sigma(11,Q106);\sigma(01,Q107);\sigma(11,Q108);\sigma(01,Q109);%
\sigma(10,Q1010);\sigma(10,Q1011);\sigma(00,Q1012);\sigma(10,Q1013);\sigma(00,Q1014);}{}; 
}; 
\path (47,0) coordinate (start);
\sigmatwo(10,start); \draw[gray](start)--(59,12);
\end{tikzpicture}
\caption{Direct product of $\theta:\;0\mapsto 0\,1\,0,\; 1\mapsto 1\,1\,0\,1\,0$, and second iterate of $\binom{1}{0}$.}
\end{figure}

\section{ Balancing with weights}\label{sec: weights}

To tackle the problem posed at the end of the previous section we generalize the notion of balanced blocks.

A \emph{weight vector} $w^\T=(w_0, w_1)$  consists of two positive reals $w_0, w_1$.
A block $\hA=\AB$ consisting of two words $A$ and $B$ is called a \emph{$w$-balanced block} if
$$ N_0(A)w_0+N_1(A)w_1=N_0(B)w_0+N_1(B)w_1,$$
which we can shortly write as $N(A)^\T w=N(B)^\T w$.

Observe that this generalizes the notion of balancedness in the sense that
 $\hA$ balanced implies $\hA$ is $w$-balanced, irrespective of the choice of $w$. Throughout the paper the eigenvalues of $M_\tht$ will be denoted by $\lambda^o$ and $\lambda$, where $\lambda$ is the eigenvalue with the largest modulus, also called PF eigenvalue.

\begin{lemma}\label{lem:w-PF}
\noindent Let $\theta$ be a primitive substitution with matrix $M_\tht$.
Let $w$ be a Perron Frobenius right eigenvector of $M_\tht$: $M_\tht w=\lambda w$.
 Then $\hA$ $w$-balanced $\Leftrightarrow$ $\hat{\tht}(\hA)$ $w$-balanced.
\end{lemma}

\noindent \emph{Proof:} Recall that $M_\tht(i,j)$ equals the number of times the letter $j$ occurs in $\theta(i)$. Thus
 $$N(\tht(A))^\T\,w= N(A)^\T\,M_\tht\,w= \lambda\,N(A)^\T w= \lambda\,N(B)^\T w= N(B)^\T\,M_\tht\,w=N(\tht(B))^\T\,w.\qed$$

It will be useful for us to have an explicit expression for the Perron Frobenius right eigenvector.
We denote by $|w|$ the length of a word $w$.

\begin{lemma}\label{lem:shape-PF}
\noindent Let $\theta$ be a primitive substitution with matrix $M_\tht$.
Let $\lambda^o,\lambda$ with $|\lambda^o|<\lambda$ be the eigenvalues of  $M_\tht$.
 Then $w^\T=(|\theta(0)|-\lambda^o,|\theta(1)|-\lambda^o)$ is a right eigenvector for $\lambda$.
 \end{lemma}

\noindent \emph{Proof:} By Cayley Hamilton we obtain $(M_\tht-\lambda\,I\!d)(M_\tht-\lambda^o\,I\!d)\oo=\zz.\qed$

\medskip

We continue with the substitution $\theta:\;0\mapsto 0\,1\,0\quad 1\mapsto 1\,1\,0\,1\,0$; here $M_\tht=\Mth$.
A right eigenvector $w$ is $w^\T=(1,2)^\T$.
So  the `double' sequence $\Delta^*=\uv$ starts with the irreducible $w$-balanced block
\begin{center}
\begin{tikzpicture}[scale=1.1,rounded corners=2]
\path (0.25,0) coordinate (HS);   \path (0,0.4) coordinate (VS); \path (VS)++(VS) coordinate (VS2); \path (HS)++(HS) coordinate (HS2);
\path (VS)++(HS2) coordinate  (PVH2);   \path (PVH2)++(VS) coordinate  (PVHV2);
  \path (PVH2)++(HS) coordinate  (PHVH2);   \path (PHVH2)++(VS) coordinate  (PHVHV2);
\draw (0,0)--(VS2);  \draw (HS2)--(PVH2)--(PHVH2)--(PHVHV2);
\path (0.125,0.2) coordinate  (PL);   \path (PL)++(VS) coordinate  (PLV);
\node at (PL) { $\;\;1\, 1$};  \node at (PLV) {$\;\;\;\;\,0\, 1\, 0$}; \node at (-0.6,0.5) {$\hC=$};
\node at (3,0.4) {directly followed by\;};
\path (6,0) coordinate (startB);
\path (0.25,0) coordinate (HS);   \path (0,0.4) coordinate (VS); \path (VS)++(VS) coordinate (VS2);
\path (HS)++(HS) coordinate (HS2); \path (HS2)++(HS) coordinate (HS3);
\path (startB) coordinate  (P0);   \path (P0)++(VS) coordinate (P1); \path (P1)++(HS) coordinate (P2); \path (P2)++(VS) coordinate (P3);
\draw (startB)++(HS3)--++(VS2);  \draw (P0)--(P1)--(P2)--(P3);
\path (startB)++(0.125,0.2) coordinate  (PL);   \path (PL)++(VS)++(HS) coordinate  (PLV); \path (startB)++(-0.6,0.5) coordinate (Text);
\path (startB)++(0.9,0.4) coordinate (end);
\node at (PL) {$\;\;\;\;\,0\, 1\, 0$ };  \node at (PLV) {$\;\;1\, 1$}; \node at (Text) {$\hD=$}; \node at (end) {$.$};
\end{tikzpicture}
\end{center}
By Lemma \ref{lem:w-PF}, we can decompose the $\hat{\theta}$-image of these two blocks in irreducible $w$-balanced blocks:
\begin{equation*}
\hat{\theta}(\hC) =
\begin{tikzpicture}[scale=1,rounded corners=2.9]
\def\Ahat(#1){
\path (#1) coordinate (P0);
\path (0.05,0) coordinate (XH);  \path (0,0.05) coordinate (YH);
\path (0.3,0) coordinate (XMH);  \path (0,0.3) coordinate (YMH);
\path (XMH)++(XH) coordinate (X); \path (YMH)++(YH) coordinate (Y);  \path (YMH)++(YH)++(YH) coordinate (YPH);
\path (P0)++(XMH)++(X) coordinate (P2);  \path (P0)++(Y)++(Y) coordinate (P1);
\path (P2)++(YPH) coordinate (P3); \path (P3)++(X) coordinate (P4); \path (P4)++(YMH) coordinate (P5);
\draw (P0)--(P1);   \draw (P2)--(P3)--(P4)--(P5);
\path (P0)++(0.125,0.2) coordinate (b1); \path (b1)++(X) coordinate (b2); \node at (b1) { $1$}; \node at (b2) {$1$};
\path (b1)++(Y) coordinate  (a1);   \path (b2)++(Y) coordinate  (a2); \path (a2)++(X) coordinate  (a3);
\node at (a1) { $0$};  \node at (a2) {$1$}; \node at (a3) { $0$}; } 
\def\Bhat(#1){
\path (0.05,0) coordinate (XH); \path(0,0.05) coordinate (YH); \path(0,0.05) coordinate (extra);
\path (#1)++(XH) coordinate (Q0); \path (0.3,0)++(XH) coordinate (X); \path (0,0.3) coordinate (YMH);
\path (YMH)++(YH) coordinate (Y); \path (Q0)++(YMH)++(extra) coordinate (Q1);  \path (YMH)++(YH)++(YH) coordinate (YPH);
\path (Q1)++(X) coordinate (Q2);  \path (Q2)++(Y) coordinate (Q3);  
\path (#1)++(X)++(X)++(X) coordinate (Q4); \path (Q4)++(Y)++(Y) coordinate (Q5);
\draw (Q0)--(Q1)--(Q2)--(Q3);   \draw (Q4)--(Q5);
\path (Q0)++(0.125,0.2) coordinate  (b1);   \path (b1)++(X) coordinate  (b2); \path (b2)++(X) coordinate  (b3);
\node at (b1) { $0$};  \node at (b2) {$1$}; \node at (b3) {$0$};
\path (b2)++(Y) coordinate  (a1);   \path (b3)++(Y) coordinate  (a2);
\node at (a1) { $1$};  \node at (a2) {$1$}; } 
\def\zenhat(#1,#2){
\path (0.05,0) coordinate (XH);  \path(0,0.05) coordinate (YH); \path(0,0.05) coordinate (extra);
\path (#1)++(XH) coordinate (Q0); \path (0.3,0)++(XH) coordinate (X); \path (0,0.3) coordinate (YMH);
\path (YMH)++(YH) coordinate (Y); \path (Q0)++(YMH)++(extra) coordinate (Q1); \path (YMH)++(YH)++(YH) coordinate (YPH);
\path (Q1)++(X) coordinate (Q2); \path (Q2)++(Y) coordinate (Q3);  
\path (#1) coordinate (P0); \path (P0)++(XMH)++(XH)++(XH) coordinate (P2); \path (P2)++(YPH) coordinate (P3);
\path (P3)++(X) coordinate (P4); \path (P4)++(YMH) coordinate (P5);
\draw (Q0)--(Q1)--(Q2)--(Q3);    \draw (P2)--(P3)--(P4)--(P5);
\path (Q0)++(0.125,0.2)++(0.03,0) coordinate  (b1);   \node at (b1) { $#2$};  
\path (b1)++(X)++(Y)++(0.04,0) coordinate  (a1);   \node at (a1) { $#2$};} 
\path (0,0) coordinate (L1); \path (L1)++(0.7,0) coordinate (L2); \path (L2)++(1.15,0) coordinate (L3); \path (L3)++(0.7,0) coordinate (L4);
\path (L4)++(.45,0) coordinate (L5); \path (L5)++(0.45,0) coordinate (L6);
\Ahat(L1);  \Bhat(L2); \Ahat(L3); \zenhat(L4,0); \zenhat(L5,1); \zenhat(L6,0);
\end{tikzpicture}\!\!\!\!,\qquad
\hat{\theta}(\hD) =
\begin{tikzpicture}[scale=1,rounded corners=2.9]
\def\Ahat(#1){
\path (#1) coordinate (P0);
\path (0.05,0) coordinate (XH);  \path (0,0.05) coordinate (YH);
\path (0.3,0) coordinate (XMH);  \path (0,0.3) coordinate (YMH);
\path (XMH)++(XH) coordinate (X); \path (YMH)++(YH) coordinate (Y);  \path (YMH)++(YH)++(YH) coordinate (YPH);
\path (P0)++(XMH)++(X) coordinate (P2);  \path (P0)++(Y)++(Y) coordinate (P1);
\path (P2)++(YPH) coordinate (P3); \path (P3)++(X) coordinate (P4); \path (P4)++(YMH) coordinate (P5);
\draw (P0)--(P1);   \draw (P2)--(P3)--(P4)--(P5);
\path (P0)++(0.125,0.2) coordinate (b1); \path (b1)++(X) coordinate (b2); \node at (b1) { $1$}; \node at (b2) {$1$};
\path (b1)++(Y) coordinate  (a1);   \path (b2)++(Y) coordinate  (a2); \path (a2)++(X) coordinate  (a3);
\node at (a1) { $0$};  \node at (a2) {$1$}; \node at (a3) { $0$}; } 
\def\Bhat(#1){
\path (0.05,0) coordinate (XH); \path(0,0.05) coordinate (YH); \path(0,0.05) coordinate (extra);
\path (#1)++(XH) coordinate (Q0); \path (0.3,0)++(XH) coordinate (X); \path (0,0.3) coordinate (YMH);
\path (YMH)++(YH) coordinate (Y); \path (Q0)++(YMH)++(extra) coordinate (Q1);  \path (YMH)++(YH)++(YH) coordinate (YPH);
\path (Q1)++(X) coordinate (Q2);  \path (Q2)++(Y) coordinate (Q3);  
\path (#1)++(X)++(X)++(X) coordinate (Q4); \path (Q4)++(Y)++(Y) coordinate (Q5);
\draw (Q0)--(Q1)--(Q2)--(Q3);   \draw (Q4)--(Q5);
\path (Q0)++(0.125,0.2) coordinate  (b1);   \path (b1)++(X) coordinate  (b2); \path (b2)++(X) coordinate  (b3);
\node at (b1) { $0$};  \node at (b2) {$1$}; \node at (b3) {$0$};
\path (b2)++(Y) coordinate  (a1);   \path (b3)++(Y) coordinate  (a2);
\node at (a1) { $1$};  \node at (a2) {$1$}; } 
\def\zenhat(#1,#2){
\path (#1) coordinate (Q0); \path (0.35,0) coordinate (X); \path (0,0.35) coordinate (Y);
\path (Q0)++(Y)++(Y) coordinate (Q1);
\path (#1) coordinate (P0); \path (P0)++(X) coordinate (P2); \path (P2)++(Y)++(Y) coordinate (P3);
\draw (Q0)--(Q1);    \draw (P2)--(P3);
\path (Q0)++(0.125,0.2)++(0.03,0) coordinate (b1);   \node at (b1) { $#2$};
\path (b1)++(Y) coordinate (a1);   \node at (a1) { $#2$};} 
\path (0,0) coordinate (L1); \path (L1)++(1.15,0) coordinate (L2); \path (L2)++(0.7,0) coordinate (L3); \path (L3)++(1.15,0) coordinate (L4);
\path (L4)++(.45,0) coordinate (L5); \path (L5)++(0.45,0) coordinate (L6);
 \Bhat(L1);  \Ahat(L2); \Bhat(L3); \zenhat(L4,0); \zenhat(L5,1); \zenhat(L6,0);
\end{tikzpicture}.\vspace*{-.2cm}
\end{equation*}
\noindent Note the occurrence of the `skewed' $w$-balanced blocks $\hat{0}_{\natural}=:{\footnotesize \skewzero}$ and $\hat{1}_{\natural}=:{\footnotesize \skewone}$.

We see that  $\Delta^*$ is generated by the substitution $\hat{\theta}$ given by:
\begin{equation*}
\begin{split}\hat{\theta}(\hat{0})=\hat{0}\,\hat{1}\,\hat{0},\quad\hat{\theta}(\hat{1})=\hat{1}\,\hat{1}\,\hat{0}\,\hat{1}\,\hat{0},
\quad\hat{\theta}(\hat{0}_{\natural})=\hat{0}_{\natural}\,\hat{1}_{\natural}\,\hat{0}_{\natural},
\quad\hat{\theta}(\hat{1}_{\natural})=\hat{1}_{\natural}\,\hat{1}_{\natural}\,\hat{0}_{\natural}\,\hat{1}_{\natural}\,\hat{0}_{\natural},\\
\quad\hat{\theta}(\hC)=\hC\,\hD\,\hC\,\hat{0}_{\natural}\,\hat{1}_{\natural}\,\hat{0}_{\natural},
\quad\hat{\theta}(\hD)=\hD\,\hC\,\hD\,\hat{0}\,\hat{1}\,\hat{0}.\hspace*{3.7cm}
\end{split}
\end{equation*}
However,... we tacitly assumed that the blocks keep their `shape'!
Actually, if $A=a_1\dots a_m$, and $B=b_1\dots b_n$, then a $w$-balanced block has a shape as

\medskip

\begin{tikzpicture}[scale=1,rounded corners=2.9]
\path (0,0) coordinate (Q0);   \path (Q0)++(0,0.41) coordinate (Q1);
\path (Q1)++(1.2,0) coordinate (Q2);  \path (Q2)++(0,0.41) coordinate (Q3);
\draw (Q0)--(Q1)--(Q2)--(Q3);
\path (6,0) coordinate (R0);   \path (R0)++(0,0.41) coordinate (R1);
\path (3.1,0)++(0,0.41) coordinate (R2); \path (R2)++(0,0.41) coordinate (R3);
\draw (R0)--(R1)--(R2)--(R3);
\path (0.125,0.2) coordinate (b1);
\path (b1)++(0,0.41) coordinate (a1);
  \node at (a1) {\hspace*{4cm}$a_1\,.\,.\,.\,.\,.\,.a_m$};
  \node at (b1) {\hspace*{5.8cm}$b_1\dots b_sb_{s+1}\dots b_{s+m}b_{s+m+1}\dots b_{n-1}b_n$};
  \node at (6.2, 0.2) {.};
\end{tikzpicture}

\noindent So a $w$-balanced block has two more parameters, a \emph{left slope} $s$ and a \emph{right slope} $m-n+s$.

\smallskip

\emph{By definition} the left slope of $\hat{\theta}(\hA)$ is equal to the left slope of $\hA$.

\begin{lemma}\label{lem:slope}
\noindent Let $\theta$ be a primitive substitution of Salem type. Let $\hA=\AB$ be a $w$-balanced block. Then the right slope of $\hat{\theta}(\AB)$ is equal to the right slope of $\AB$.\ts
\end{lemma}

\noindent \emph{Proof:} Suppose that $\lambda^o=1$, and let $s$ be the left slope of $\hA$.
We have to show that $|\theta(A)|-|\theta(B)|+s=|A|-|B|+s$, which holds if and only if $|\theta(A)|-|A|=|\theta(B)|-|B|.$
The left hand side equals
\begin{equation*}
\begin{split}
|\theta(A)|-|A|& =N_0(A)|\theta(0)|+N_1(A)|\theta(1)|-N_0(A)-N_1(A)\\
               & =N_0(A)(|\theta(0)|-1)+N_1(A)(|\theta(1)|-1) = N(A)^\T w,
\end{split}
\end{equation*}
where  $w$ is a Perron Frobenius right eigenvector of $M_\tht$ (using Lemma \ref{lem:shape-PF}).
Since $\AB$ is balanced  $N(A)^\T w= N(B)^\T w$, and it follows that $|\theta(A)|-|A|=|\theta(B)|-|B|.$

In case  $\lambda^o=-1$, one replaces $\tht$ by  $\tht^2$, and proceeds as above.\qed

\section{Computing the coincidence density}\label{sec: comp}

The following general lemma will be useful.

\begin{lemma}\label{lem:subseq}
\noindent Let $(a_n)$ be a sequence of real numbers such that $\lim_{n\rightarrow\infty}\frac1{n}\sum_{i=1}^na_i=a$ exists.
Let $(n_k)$ and $(t_k)$ be two sequences with $t_k\ge 0,\; n_k\rightarrow\infty$ and $\limsup n_k/t_k<\infty$. Then
$$\lim_{k\rightarrow\infty}\frac1{t_k}\sum_{i=n_k+1}^{n_k+t_k}a_i=a.$$
\end{lemma}

\noindent \emph{Proof:} Let $\eps>0$. For all large $n$ we will have $n(a-\eps)\le\sum_{i=1}^{n}a_i  \le n(a+\eps).$
So
\begin{equation*}
(n_k+t_k)(a-\eps)-n_k(a+\eps)\le\sum_{i=1}^{n_k+t_k}a_i-\sum_{i=1}^{n_k}a_i \le (n_k+t_k)(a+\eps)-n_k(a-\eps),
\end{equation*}
which is equivalent to\\[-.4cm]
\begin{equation*}
 t_ka-\eps(t_k+2n_k) \le\sum_{i=n_k+1}^{n_k+t_k}a_i\le t_ka-\eps(t_k+2n_k).
\end{equation*}
Dividing through by $t_k$, and letting $k\rightarrow\infty$, the result follows by next letting $\eps\downarrow 0$.\qed\\

 Our main result is

\begin{theorem}\label{th:main}
\noindent Let $\theta$ be a primitive substitution on $\{0,1\}^*$ of Salem type with $\theta(0)=0\dots$, and $\theta(1)=1\dots$. Let $w$
be a  Perron-Frobenius right eigenvector of $M_\tht$. Suppose the set of $w$-balanced blocks is finite, and contains both $\hat{0}$ and $\hat{0}_{\natural}$. Then the coincidence density of $\theta$ does not exist.
\end{theorem}

\noindent \emph{Proof:} Suppose that $\delta_c$ \emph{does} exist. Let $j_0$ be a position in the double sequence $\hat{x}$ where
$\hat{0}$ does occur. Let $A=x_1\dots x_{j_0-1}$. It follows from Lemma~\ref{lem:slope} that
for all $k\ge 1$ the block $\hat{\tht}^k(\hat{0})$ does occur in $\hat{x}$ at position $n_k:=|\tht^k(A)|$. Let $t_k:=|\tht^k(0)|$. From
the Perron-Frobenius theorem one deduces that $\limsup n_k/t_k<\infty$. We can therefore apply  Lemma~\ref{lem:subseq}, and
since we have $t_k$ coincidences in the block $\hat{\tht}^k(\hat{0})$, it follows that $\delta_c=1$. Similarly, let $j_1$ be a position
where  $\hat{0}_{\natural}$ occurs, and take $n_k$ (with $j_0$ replaced by $j_1$) as before. Ignoring the lower left and the upper right letter of $\hat{\tht}^k(\hat{0}_{\natural})$, take  $t_k=|\theta^k(0)|-2$. It remains to determine the number of coincidences in the block $\hat{\tht}^k(\hat{0}_{\natural})$ (with the two letters removed). For any $k$ this block contains a positive fraction of blocks $\hat{\tht}(\hat{0}_{\natural})$, so not
to have a contradiction with $\delta_c=1$, no \emph{non}-coincidences are allowed for $\hat{\tht}(\hat{0}_{\natural})$. This requirement
leads to $\tht(0)=0\,0\dots 0$, which contradicts primitivity of $\theta$. Conclusion:  $\delta_c$ does \emph{not} exist. \qed

\smallskip

Applying this result to  the substitution $\theta:\;0\mapsto 0\,1\,0,\, 1\mapsto 1\,1\,0\,1\,0$ from the previous section,  we find that the coincidence
density does not exists here, nor will the frequencies of letters on the diagonal $\Delta^*(0,1)$ exist.


\section{Self-similarity can increase regularity}\label{sec: selfsim}

 As in \cite{Solom} one obtains a \emph{self similar} 1D tiling from a substitution $\tht$ by associating to each letter $i$ an interval of  length $w_i$, where $w$ is a right Perron-Frobenius eigenvector of the substitution matrix. This induces a 2D self-similar tiling from   the direct product $\thtt$, by associating to each pair of letters $(i,j)$ the  $w_i\times w_j$ rectangle. Let us reconsider our  example $\theta:\;0\mapsto 0\,1\,0,\; 1\mapsto 1\,1\,0\,1\,0$, choosing $w^\T=(1,2)^\T$.
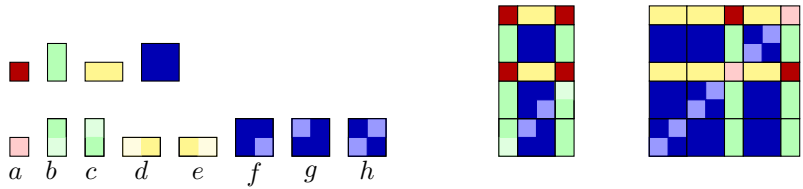
\begin{figure}[h!]
\centering
\begin{tikzpicture}[scale=.25,rounded corners=0]
\def\brickzerozero(#1){\path (#1) coordinate (P0);  \draw [fill=Color00, draw=black]  (P0) rectangle +(1,1)}
\def\brickzeroone(#1){\path (#1) coordinate (P0);  \draw [fill=Color01, draw=black]  (P0) rectangle +(1,2)}
\def\brickonezero(#1){\path (#1) coordinate (P0); \draw [fill=Color10, draw=black]  (P0) rectangle +(2,1)}
\def\brickoneone(#1){\path (#1) coordinate (P0); \draw [fill=Color11, draw=black]  (P0) rectangle +(2,2)}
\def\brickzerozerod(#1){\path (#1) coordinate (P0);  \draw [fill=Color00d, draw=black]  (P0) rectangle +(1,1)}
\def\brickzerooned(#1){\path (#1) coordinate (P0);  \draw [fill=Color01d, draw=black]  (P0) rectangle +(1,2)}
\def\brickoneoned(#1){\path (#1) coordinate (P0); \draw [fill=Color11d, draw=black]  (P0) rectangle +(2,2)}
\def\brickzeroonedd(#1){\path (#1) coordinate (P0);  \draw [fill=Color01]  (P0) rectangle +(1,2);
    \draw [fill=Color01d,draw=Color01] (P0) rectangle +(1,1); \draw [draw=black]  (P0) rectangle +(1,2) }
\def\brickzeroonedu(#1){\path (#1) coordinate (P0);  \path (P0)++(0,1) coordinate (P1); \draw [fill=Color01]  (P0) rectangle +(1,2);
    \draw [fill=Color01d,draw=Color01] (P1) rectangle +(1,1); \draw [draw=black]  (P0) rectangle +(1,2) }
\def\brickonezerodl(#1){\path (#1) coordinate (P0); \draw [fill=Color10]  (P0) rectangle +(2,1);
    \draw [fill=Color10d, draw=Color10] (P0) rectangle +(1,1); \draw [draw=black]  (P0) rectangle +(2,1) }
\def\brickonezerodr(#1){\path (#1) coordinate (P0); \path (P0)++(1,0) coordinate (P1); \draw [fill=Color10]  (P0) rectangle +(2,1);
    \draw [fill=Color10d, draw=Color10] (P1) rectangle +(1,1); \draw [draw=black]  (P0) rectangle +(2,1) }
\def\brickoneonedu(#1){\path (#1) coordinate (P0); \path (P0)++(0,1) coordinate (P1); \draw [fill=Color11]  (P0) rectangle +(2,2);
    \draw [fill=Color11d, draw=Color11] (P1) rectangle +(1,1); \draw [draw=black]  (P0) rectangle +(2,2) }
\def\brickoneonedd(#1){\path (#1) coordinate (P0); \path (P0)++(1,0) coordinate (P1); \draw [fill=Color11]  (P0) rectangle +(2,2);
     \draw [fill=Color11d, draw=Color11] (P1) rectangle +(1,1); \draw [draw=black]  (P0) rectangle +(2,2) }
\def\brickoneonedt(#1){\path (#1) coordinate (P0); \path (P0)++(1,1) coordinate (P1); \draw [fill=Color11]  (P0) rectangle +(2,2);
     \draw [fill=Color11d, draw=Color11] (P0) rectangle +(1,1);\draw [fill=Color11d, draw=Color11] (P1) rectangle +(1,1);
     \draw [draw=black]  (P0) rectangle +(2,2) }
\path (0,4) coordinate (start); \path (2,4) coordinate (next1); \path (4,4) coordinate (next2);  \path (7,4) coordinate (next3);
\brickzerozero(start); \brickzeroone(next1); \brickonezero(next2); \brickoneone(next3);
\path (0,0) coordinate (start00);  \path (2,0) coordinate (next01dd);  \path (4,0) coordinate (next01du); \path (6,0) coordinate (next10dl);
\path (9,0) coordinate (next10dr); \path (15,0) coordinate (next11dd); \path (12,0) coordinate (next11du);\path (18,0) coordinate (next11dt);
\brickzerozerod(start00); \brickzeroonedd(next01dd); \brickzeroonedu(next01du); \brickonezerodl(next10dl); \brickonezerodr(next10dr);
\brickoneonedu(next11dd); \brickoneonedd(next11du);  \brickoneonedt(next11dt);
\node at (0.3,-.9) {$a$ };  \node at (2.2,-.75) {$b$ };  \node at (4.3,-.9) {$c$ };  \node at (7.0,-.75) {$d$ }; \node at (10.0,-.9) {$e$ };
\node at (13.0,-.95) {$f$ };\node at (16.0,-.95) {$g$ };\node at (19.0,-.75) {$h$ };
\def\sigma(#1,#2){\path (#2) coordinate (P0);
\path (P0)++(0,0) coordinate  (P010); \path (P0)++(1,0) coordinate  (P011); \path (P0)++(3,0) coordinate  (P012);
\path (P0)++(0,2) coordinate  (P013); \path (P0)++(1,2) coordinate  (P014); \path (P0)++(3,2) coordinate  (P015);
\path (P0)++(0,4) coordinate  (P016); \path (P0)++(1,4) coordinate  (P017); \path (P0)++(3,4) coordinate  (P018);
\path (P0)++(0,5) coordinate  (P019); \path (P0)++(1,5) coordinate  (P0110); \path (P0)++(3,5) coordinate  (P0111);
\path (P0)++(0,7) coordinate  (P0112); \path (P0)++(1,7) coordinate  (P0113); \path (P0)++(3,7) coordinate  (P0114);
\ifthenelse{#1=0170}
{\brickzeroonedd(P010);\brickoneonedu(P011);\brickzeroone(P012); \brickzeroone(P013);\brickoneonedd(P014);\brickzeroonedu(P015);%
\brickzerozero(P016);\brickonezero(P017);\brickzerozero(P018); \brickzeroone(P019);\brickoneone(P0110);\brickzeroone(P0111);%
\brickzerozero(P0112);\brickonezero(P0113);\brickzerozero(P0114);}{};                                     
\path (P0)++(0,0) coordinate (P110); \path (P0)++(2,0) coordinate (P111); \path (P0)++(4,0) coordinate (P112);
\path (P0)++(5,0) coordinate (P113); \path (P0)++(7,0) coordinate (P114); \path (P0)++(0,2) coordinate (P115);
\path (P0)++(2,2) coordinate (P116); \path (P0)++(4,2) coordinate (P117); \path (P0)++(5,2) coordinate (P118);
\path (P0)++(7,2) coordinate (P119); \path (P0)++(0,4) coordinate (P1110); \path (P0)++(2,4) coordinate (P1111);
\path (P0)++(4,4) coordinate (P1112); \path (P0)++(5,4) coordinate (P1113); \path (P0)++(7,4) coordinate (P1114);
\path (P0)++(0,5) coordinate (P1115); \path (P0)++(2,5) coordinate (P1116); \path (P0)++(4,5) coordinate (P1117);
\path (P0)++(5,5) coordinate (P1118);\path (P0)++(7,5) coordinate (P1119);  \path (P0)++(0,7) coordinate (P1120);
\path (P0)++(2,7) coordinate (P1121);\path (P0)++(4,7) coordinate (P1122); \path (P0)++(5,7) coordinate (P1123);\path (P0)++(7,7) coordinate (P1124);
\ifthenelse{#1=1172}
{\brickoneonedt(P110);\brickoneone(P111);\brickzeroone(P112);\brickoneone(P113);\brickzeroone(P114);%
\brickoneone(P115);\brickoneonedt(P116);\brickzeroone(P117);\brickoneone(P118);\brickzeroone(P119);%
\brickonezero(P1110);\brickonezero(P1111);\brickzerozerod(P1112);\brickonezero(P1113);\brickzerozero(P1114);%
\brickoneone(P1115);\brickoneone(P1116);\brickzeroone(P1117);\brickoneonedt(P1118);\brickzeroone(P1119);%
\brickonezero(P1120);\brickonezero(P1121);\brickzerozero(P1122);\brickonezero(P1123);\brickzerozerod(P1124);}{}; 
}; 
\path (26,0) coordinate (start); \path (34,0) coordinate (next); \sigma(0170,start); \sigma(1172,next);
\end{tikzpicture}
\caption{Four tiles, 8 diagonal tiles, and the first iterate of the $b$-tile and the $h$-tile.}\label{fig:selfsim}
\end{figure}
The diagonal of the 2D self-similar tiling associated to this example can be generated by a substitution $\delta$ on eight letters, corresponding to the ways in which the diagonal intersects a tile: $a$ for the $1\times 1$ tile,  $b$ and $c$ for the $1\times 2$ tile, etc., see Figure~\ref{fig:selfsim}.
The substitution $\dd$ is  given by
$$\dd(a)=\dd(c)=\dd(e)=a\,h\,a,\; \dd(b)=\dd(f)=b\,g\,f\,c,\; \dd(d)=\dd(g)=d\,f\,g\,e,\; \dd(h)=h\,h\,a\,h\,a.$$
We can split the word $\dd(h)$ into the words $h\,h$ and $a\,h\,a$, and the letter $h$ into two letters $h=\ha\hb$, obtaining a new substitution $\ddp$ on nine letters generating the diagonal, given by
$$\ddp(a)=\ddp(c)=\ddp(e)=\ddp(\hb)=a\,\ha\,\hb\,a,\; \ddp(b)=\ddp(f)=b\,g\,f\,c,\; \ddp(d)=\ddp(g)=d\,f\,g\,e,\;
\ddp(\ha)=\ha\,\hb\,\ha\,\hb.$$
Now note that $\ddp$ is a substitution of constant length 4, and discrete. This implies that the frequency of letters in $\Delta^*(0,1)=\Delta^*(b)$
exists (the letters $a$ and $h$ will have frequencies $2/3$ and $1/3$, the others frequency 0). This phenomenon is a special case of the following general result.
\begin{proposition}\label{prop:selfsim}
\noindent Let $\tht$ be a primitive substitution on $\{0,1\}^*$ with $\theta(0)=0\dots,\; \theta(1)=1\dots$, for which $M_\tht$ has rational eigenvalues. Then the frequencies of letters on the diagonal of the 2D self-similar tiling generated by $\thtt$ do exist.
\end{proposition}

\noindent \emph{Proof:} Let $\lambda$ be the PF eigenvalue of $\tht$. Since sum and product of the eigenvalues are integers, it follows that $\lambda$ is an integer, and thus  there is an integer eigenvector, say $(p,q)^\T$, with $p,q>0$. Ignoring a scaling, we can assume that the 1D self-similar tiling generated by $\tht$ is built out of two tiles with length $p$, respectively $q$. The strategy is to cut up these two tiles in intervals of length 1, so that we may describe the diagonal of the direct product by the corresponding $1\times 1$ squares (cf.~the example above). Let $a_k,\; k=1,\dots, p$ and $b_\ell,\; \ell=1,\dots,q$ be the letters of an alphabet $\dot{A}$ of size $p+q$. Define the  map $\gamma:\{0,1\}^*\rightarrow\dot{A}^*$ by
$$\gamma(0)=a_1a_2\dots a_p,\quad \gamma(1)=b_1b_2\dots b_q.$$
Next, define the substitution $\dot{\tht}:\dot{A}^*\rightarrow\dot{A}^*$ by
$$\dot{\tht}(a_k)=\gamma(\tht(0))[(k-1)\lambda+1,\,k\lambda],\; k=1,\dots p,\quad \dot{\tht}(b_\ell)=\gamma(\tht(1))[(\ell-1)\lambda+1,\,\ell\lambda],\; \ell=1,\dots,q.$$
Here we use the notation $w[m,n]:=w_m\dots w_n$, and we remark that such a definition makes sense because of the eigenvalue equation
$M_\tht \pq=\lambda\pq$. Now note that the diagonal of the direct product self-similar tiling is generated by $\dot{\tht}\times\dot{\tht}$, where $\dot{\tht}$ is a  \emph{constant length} substitution. As remarked in Section~\ref{sec: intro}, this implies that the frequencies of letters exist.\qed

\section{All eigenvalues outside the unit circle}\label{sec: outside}

 The case where the modulus of both eigenvalues of $M_\tht$ is larger than 1 is more difficult. It is conjectured in \cite{PierreM} that in this case the coincidence density will not exist.

 \smallskip

 First we discuss the case where the eigenvalues of $M_\tht$ are rational. As in the proof of Proposition~\ref{prop:selfsim} this implies that the smallest eigenvalue $\lambda^o$ is an integer. Let $w$ be a right PF vector of $M_\tht$. The set of $w$-balanced blocks will be either empty or infinite. This follows from an extension of Lemma~\ref{lem:slope}: the left and the right slope of
$w$-balanced blocks are multiplied by $\lambda^o$ when passing from $\hA$ to $\hat{\theta}(\hA)$. Thus the technique of $w$-balanced blocks is no longer useful. However, for the \emph{self-similar} 2D tilings generated by $\tht$ Proposition~\ref{prop:selfsim} applies.

\smallskip

Secondly we consider irrational eigenvalues.  With Lemma~\ref{lem:shape-PF} one deduces from this irrationality that the only $w$-balanced blocks (with the PF weight vector $w$) are the balanced blocks. Using the condition $|\lambda^o|>1$, one can show that the set of balanced blocks
is either empty or infinite.

As an example of a substitution with irrational eigenvalues we consider
$$  \tht:\; 0\mapsto 0\,0\,0\,0\,1, \quad 1\mapsto 1\,1\,1\,0,$$
which has eigenvalues $\frac72-\frac12\sqrt{5}$ and $\frac72+\frac12\sqrt{5}$.

\begin{proposition}\label{prop:nobalance}
\noindent The set of balanced blocks of the substitution $\tht$ given by $0\mapsto 0\,0\,0\,0\,1, \, 1\mapsto 1\,1\,1\,0$ is empty.
\end{proposition}

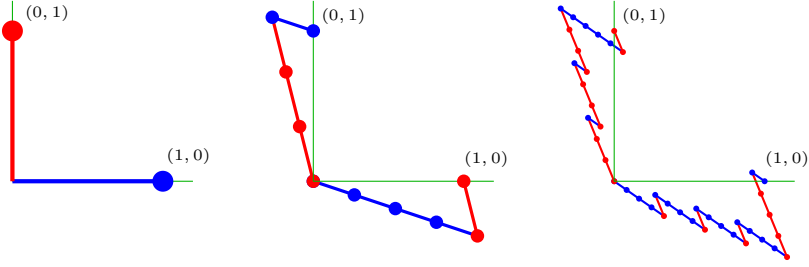
\begin{figure}[t!]
\centering
\begin{tikzpicture}[scale=2]
\path (3*\EL,-\EL) coordinate (a);        \path (-\EL,4*\EL) coordinate (b);
\path (10*\ELL,-7*\ELL) coordinate (tha); \path (-7*\ELL,17*\ELL) coordinate (thb);
\def\step(#1,#2){
\ifthenelse{#1=0} {\node[fill=blue,circle, inner sep=1.8pt] at (#2)  {};\draw [blue, very thick] (#2)--++(a)}{};
\ifthenelse{#1=1} {\node[fill=red,circle, inner sep=1.8pt] at (#2)  {};\draw [red, very thick] (#2)--++(b)}{};  }
\def\stepp(#1,#2){
\ifthenelse{#1=0} {\node[fill=blue,circle, inner sep=0.8pt] at (#2)  {};\draw [blue,  thick] (#2)--++(tha)}{};
\ifthenelse{#1=1} {\node[fill=red,circle, inner sep=0.8pt] at (#2)  {};\draw [red,  thick] (#2)--++(thb)}{};  }
\def\sigma(#1,#2){
\ifthenelse{#1=0}
{\step(0,#2); \path (#2)+(a) coordinate (D); \step(0,D); \path (D)+(a) coordinate (D); \step(0,D);%
 \path (D)+(a) coordinate (D); \step(0,D); \path (D)+(a) coordinate (D); \step(1,D);}{};  
\ifthenelse{#1=1}
{ \step(1,#2); \path (#2)+(b) coordinate (D); \step(1,D); \path (D)+(b) coordinate (D); \step(1,D);%
 \path (D)+(b) coordinate (D); \step(0,D); }{};  
}; 
\def\sigmaa(#1,#2){
\ifthenelse{#1=0}
{\stepp(0,#2); \path (#2)+(tha) coordinate (D); \stepp(0,D); \path (D)+(tha) coordinate (D); \stepp(0,D);%
 \path (D)+(tha) coordinate (D); \stepp(0,D); \path (D)+(tha) coordinate (D); \stepp(1,D);}{};  
\ifthenelse{#1=1}
{ \stepp(1,#2); \path (#2)+(thb) coordinate (D); \stepp(1,D); \path (D)+(thb) coordinate (D); \stepp(1,D);%
 \path (D)+(thb) coordinate (D); \stepp(0,D); }{};  
}; 
\path (0,0) coordinate (start);
\sigma(0,start); \sigma(1,start); \draw[green!70!black](start)--++(1.2,0); \draw[green!70!black](start)--++(0,1.2);
\node[fill=blue,circle, inner sep=1.8pt] at (0,1) {}; \node[fill=red,circle, inner sep=1.8pt] at (1,0) {};
\def\sigmatwo(#1,#2){
\ifthenelse{#1=0}
{\sigmaa(0,#2); \path (#2)+(a) coordinate (E); \sigmaa(0,E); \path (E)+(a) coordinate (E); \sigmaa(0,E);%
 \path (E)+(a) coordinate (E); \sigmaa(0,E); \path (E)+(a) coordinate (E); \sigmaa(1,E);}{};  
\ifthenelse{#1=1}
{ \sigmaa(1,#2); \path (#2)+(b) coordinate (E); \sigmaa(1,E); \path (E)+(b) coordinate (E); \sigmaa(1,E);%
 \path (E)+(b) coordinate (E); \sigmaa(0,E); }{};  
}; 
\path (2,0) coordinate (start);
\sigmatwo(0,start); \sigmatwo(1,start); \draw[green!70!black](start)--++(1.2,0); \draw[green!70!black](start)--++(0,1.2);
\node[fill=blue,circle, inner sep=0.8pt] at (3,0) {}; \node[fill=red,circle, inner sep=0.8pt] at (2,1) {};
\path (-2,0) coordinate (start);
\draw[green!70!black](start)--++(1.2,0); \draw[green!70!black](start)--++(0,1.2);
\node[fill=blue,circle, inner sep=2.8pt] at (-1,0) {}; \node[fill=red,circle, inner sep=2.8pt] at (-2,1) {};
\draw[blue, ultra thick](start)--++(1,0); \draw[red, ultra thick](start)--++(0,1);
\node at (1.15,0.15) {\tiny $(1,0)$ }; \node at (3.15,0.15) {\tiny $(1,0)$ }; \node at (-0.83,0.17) {\tiny $(1,0)$ };
\node at (0.2,1.1) {\tiny $(0,1)$ }; \node at (2.2,1.1) {\tiny $(0,1)$ }; \node at (-1.77,1.12) {\tiny $(0,1)$ };
\end{tikzpicture}
\caption{The first three approximants to the two curves $K(0)$ and $K(1)]$ generated by  the full representation of the substitution $0\mapsto 0\,0\,0\,0\,1, \; 1\mapsto 1\,1\,1\,0$.}\label{fig:nobalance}
\end{figure}

{\emph{Proof:}} We give a geometrical proof. We consider a geometric realization of $\tht$ given by its \emph{full representation}, cf.~\cite{D}. This means that we consider the morphism $f:\{0,1\}^*\rightarrow \R^2$ given by $f(0)=(1,0)^\T,\, f(1)=(0,1)^\T$, and the linear map $L_\tht$ with matrix $M_\tht$. Then we have the representation equation
$$ f\circ\tht= L_\tht\circ f.$$
From \cite{D} we have that curves $$K_n(i):=L_\tht^{-n}K[\tht^n(i)]$$ converge exponentially fast to a curve $K(i)$ for $i=0,1$. Here for a word $w=w_1\dots w_m$ we denote by $K[w]$  the polygonal curve passing successively through the points $(0,0), f(w_1), f(w_1w_2), \dots,$ $f(w_1\dots w_m)$.

In Figure~\ref{fig:nobalance} we display $K_n(0)$ and $K_n(1)$ for $n=0,1,2$. From the exponential convergence (speed $1/11$), and the
fact that these curves (almost) live outside the positive quadrant, it is obvious that they will never intersect, and thus there is no balanced block.\qed

\medskip

\begin{figure}[h!]
\centering
   \includegraphics[width=7cm]{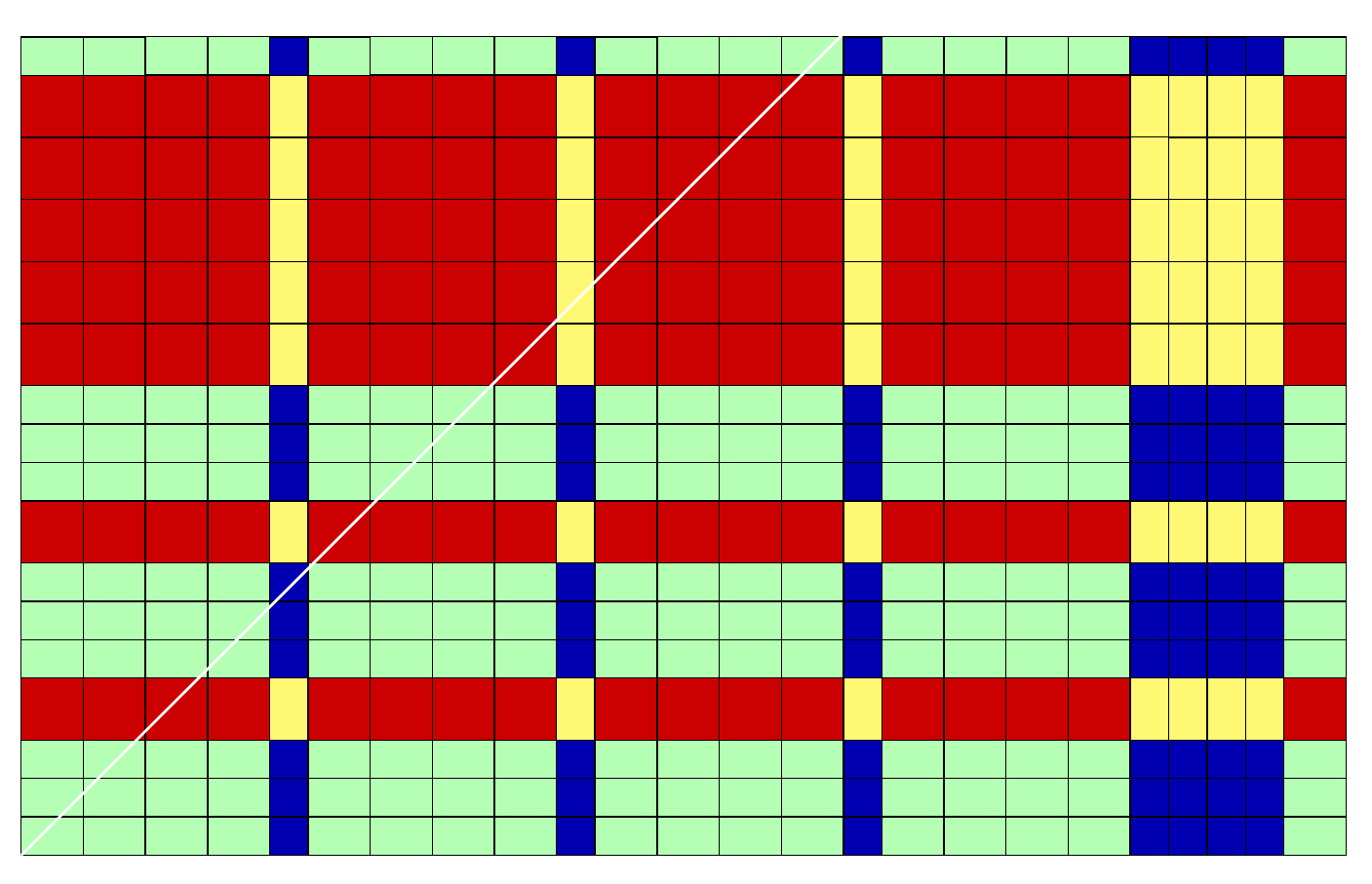}
\caption{Part of the  self-similar product tiling generated by the substitution $0\mapsto 0\,0\,0\,0\,1, \; 1\mapsto 1\,1\,1\,0$. }\label{fig:ss-nobalance}
\end{figure}

The self similar 2D-tiling associated to the substitution $0\mapsto 0\,0\,0\,0\,1, \; 1\mapsto 1\,1\,1\,0$ is illustrated with Figure~\ref{fig:ss-nobalance}. The scalings are given by the eigenvector $(\frac12+\frac12\sqrt{5},1)$. We conjecture that there are infinitely many `types' of squares and rectangles on the diagonal. At least it is clear that there are no `regeneration points' as in the Pisot case, because of Proposition~\ref{prop:nobalance}: the only $(\frac12+\frac12\sqrt{5},1)$-balanced blocks are the balanced blocks, and these do not exist.

\section{Overlap distributions}\label{sec: overlap}

The coincidence density $\delta_c$ can be considered as a special case of the so called \emph{overlap} between two sequences $x$ and $y$
defined by
$$\delta_c(x,y):=\lim_{n\rightarrow\infty}\,\frac{1}{n}{\rm Card} \{1\le k\le n: x_k=y_k\}.$$
Suppose $x$ and $y$ come from a shift invariant set $X$ of sequences, equipped with a shift invariant measure $\mu$. Then the distribution of
$\delta_c(\cdot,\cdot)$ under the measure $\mu\times\mu$ is called the \emph{overlap distribution}, introduced by Parisi (\cite{Parisi}).
It was shown in \cite{Enter-overlap} that if $\mu$ is weakly mixing, then the overlap distribution degenerates to a point mass
concentrated on $\mu([0])^2+\dots+\mu([r])^2$, where $r$ is the number of letters in the alphabet, and $[i]$ is the cylinder $\{x\in X:
x_0=i\}$. (This follows immediately from the ergodic theorem applied with $\mu\times\mu$.)

Primitive substitutions $\theta$ generate a shift
dynamical system with a unique invariant measure $\mu_\tht$. The overlap distribution of some of these is determined in \cite{Enter-overlap}
(see also \cite{Enter-fold}). For the substitution $\theta:\;0\mapsto 0\,1\,0\quad 1\mapsto 1\,1\,0\,1\,0$ studied in Section~\ref{sec: weights}, it is known that $\mu_\tht$ is weak mixing (\cite{Martin}, \cite{Host}). So for this substitution $\delta_c(x,y)$ exists for $\mu_\tht$-almost all $x$ and $y$, and for those $x$ and $y$
$$\delta_c(x,y)=\mu_\tht([0])^2+\mu_\tht([1])^2=\frac14+\frac14=\frac12.$$
By Theorem~\ref{th:main} this means that  the pair $(u^*,v^*)$ is non-generic for the overlap distribution. In fact most substitutions with matrix $\Mth$ and two fixed points will have non-generic pairs  $(u^*,v^*)$ by an application of this theorem. One exception is the substitution given by $\tht(0)=010, \tht(1)=10101$, since it generates a periodic sequence. Another is the `Kakutani' substitution (cf.~\cite{DK}) given by
$$\tht(0)=0\,0\,1,\quad \theta(1)=1\,1\,0\,0\,1.$$
Here one can show easily that $\delta_c(u^*,v^*)=\delta_c=\frac12$, so in this case the pair $(u^*,v^*)$ \emph{is} generic.

\section{ Final remarks}\label{sec: final}

We would like to point out that the phenomenon of non-existence of frequencies of symbols can already occur for substitutions of constant length. A simple example is given by
$$ 1\mapsto 1\,2\,3,\quad 2\mapsto 2\,2\,2,\quad 3\mapsto 3\,3\,3.$$

\smallskip

 Of the 12 substitutions with matrix $\Mth$ (and $\tht(0)=0\dots,\; \tht(1)=1\dots$) there  are 9 for which the coincidence density does not exist. Then there is a periodic one, and the `Kakutani' substitution mentioned in the previous section.  There is one substitution left, given by
$\theta(0)=001,\;\theta(1)=10110$. We do not know how to determine the coincidence density or frequencies of letters on the diagonal for this substitution. Nor do we know how to compute these for substitutions whose matrix has two eigenvalues outside the unit circle.


\end{document}